\documentclass[11pt,leqno]{article}

\usepackage{amsmath,amsfonts,amscd,amssymb,theorem}

\long\def\comment#1\endcomment{}

\comment
\endcomment

\makeatletter
\begingroup
\gdef\th@dotted{\normalfont\itshape
  \def\@begintheorem##1##2{%
        \item[\hskip\labelsep \theorem@headerfont ##1\ ##2.]}%
\def\@opargbegintheorem##1##2##3{%
   \item[\hskip\labelsep \theorem@headerfont ##1\ ##2\ (##3).]}}
\endgroup
\makeatother

\theoremstyle{dotted}

\newtheorem{theorem}{Theorem}[section]
\newtheorem{lemma}[theorem]{Lemma}

\newtheorem{prop}[theorem]{Proposition}
\newtheorem{corr}[theorem]{Corollary}

\makeatletter
\begingroup
\gdef\th@upshape{\normalfont
  \def\@begintheorem##1##2{%
        \item[\hskip\labelsep \theorem@headerfont ##1\ ##2.]}%
\def\@opargbegintheorem##1##2##3{%
   \item[\hskip\labelsep \theorem@headerfont ##1\ ##2\ (##3).]}}
\endgroup
\makeatother

\theoremstyle{upshape}

\newtheorem{defn}[theorem]{Definition}
\newtheorem{notation}[theorem]{Notation}

\newtheorem{exa}[theorem]{Example}

\makeatletter
\renewcommand{\subsection}{\@startsection{subsection}{2}{0pt}{-3ex
plus -1ex minus -0.2ex}{-2mm plus -0pt minus
-2pt}{\normalfont\bfseries}} 
\renewcommand{\subsubsection}{\@startsection{subsubsection}{3}{0pt}{-3ex
plus -1ex minus -0.2ex}{-2mm plus -0pt minus
-2pt}{\normalfont\bfseries}} 
\makeatother

\makeatletter
\@addtoreset{equation}{section}
\makeatother

\newcommand{\cntrct}                
{\hspace{2pt}\raisebox{1pt}{\text{$\lrcorner$}}\hspace{2pt}}

\newcommand{\proof}[1][Proof.]{\smallskip\noindent{\em #1}}
\def\endproof{\hfill\ensuremath{\square}\par\medskip}

\def\eqref#1{\thetag{\ref{#1}}}

\let\latexref=\ref
\def\ref#1{{\normalfont{\latexref{#1}}}}

\newcommand{\wt}{\widetilde}
\newcommand{\wh}{\widehat}

\newcommand{\idot}{{\:\raisebox{1pt}{\text{\circle*{1.5}}}}}
\newcommand{\hdot}{{\:\raisebox{3pt}{\text{\circle*{1.5}}}}}
\newcommand{\eps}{\varepsilon}
\renewcommand{\phi}{\varphi}

\newcommand{\vH}{\check{H}}

\newcommand{\Ext}{\operatorname{Ext}}

\newcommand{\id}{\operatorname{\sf id}}
\newcommand{\Id}{\operatorname{\sf Id}}
\newcommand{\gr}{\operatorname{\sf gr}}
\newcommand{\tr}{\operatorname{\sf tr}}
\newcommand{\colim}{\operatorname{colim}}

\newcommand{\C}{{\cal C}}

\newcommand{\hush}{\natural}

\newcommand{\Aut}{{\operatorname{Aut}}}

\newcommand{\amod}{{\text{\rm -mod}}}

\newcommand{\ppt}{{\sf pt}}

\newcommand{\lotimes}{\overset{\sf\scriptscriptstyle L}{\otimes}}

\newcommand{\Spec}{\operatorname{Spec}}

\newcommand{\Z}{{\mathbb Z}}
\newcommand{\Q}{{\mathbb Q}}
\newcommand{\N}{{\mathbb N}}
\newcommand{\F}{{\mathbb F}}
\newcommand{\R}{{\mathbb R}}

\newcommand{\calo}{\mathcal{O}}

\newcommand{\Fun}{\operatorname{Fun}}

\newcommand{\E}{\mathcal{E}}

\newcommand{\Ab}{\operatorname{Ab}}

\newcommand{\B}{\mathcal{B}}
\newcommand{\Zz}{\mathcal{Z}}

\newcommand{\W}{\mathbb{W}}

\title{Witt vectors, commutative and non-commutative}

\author{D. Kaledin\thanks{Supported by the Russian Science
    Foundation, grant 14-50-00005}}

\begin{document}

\maketitle

\tableofcontents

\section*{Introduction.}

Witt vectors occupy a curious place in the mathematical grand scheme
of things: above the radar, but only just. Most mathematicians are
vaguely aware that such a thing exists, and no more than that. Those
who do encounter Witt vectors in the course of their own work often
develop quite an attachment towards the subject, and a bit of a
personal viewpoint. Such an encounter can happen quite unexpectedly,
and from many different directions, so that these viewpoints vary
greatly. Somehow, the beautiful discovery of Witt refuses to take
its richly deserved place in history and insists on staying alive;
moreover, it even refuses to be pigeon-holed, and reveals a new side
to itself every time one takes a new look.

The goal of this paper is to give an overview of some old results on
Witt vectors and present some new ones. However, given the shifting
nature of the subject, it is perhaps prudent to drop all pretence
right away and state clearly that what is presented here is the
viewpoint of the author, as idiosyncratic and personal as anyone
else's. To get some context, we start with a bit of a historical
overview, but again, this is not history in the proper sense:
rather, a story.

The story starts in Germany in about 1936, with Oswald Teichm\"uller
and Ernst Witt. Teichm\"uller discovered a new way to think about
$p$-adic integers: instead of treating $\Z_p$ as a completion of
$\Z$ whose elements are convergent sequences modulo an equivalence
relation, he noticed that there is a completely canonical
identification between the set $\Z_p$ and the set $W(\F_p)$ of
infinite sequences $\langle a_0,a_1,\dots \rangle$ of elements in
the prime field $\F_p$. Witt then showed \cite{W} how to turn the
set $W(\F_p)$ into an abelian group and furthermore, into a
commutative ring, by constructing explicit polynomials that provide
addition and multiplication.

The key word here is ``canonical'': Teichm\"uller's representation
of $\Z_p$ was independent of any choices, thus sufficiently natural,
and it is this naturality that allows Witt universal polynomials to
exist.

However, as a byproduct of the construction, Teichm\"uller and Witt
actually obtained much more. Namely, Witt polynomials have integer
coefficients, thus make sense not only in $\F_p$ but in any
commutative ring. So, although they did not think in these terms ---
which did not even exist at the time --- what Teichm\"uler and Witt
discovered was a functor from commutative rings to commutative
rings.

For the next step, fast-forward twenty years and move to France,
where modern algebraic geometry was being created. One of the
driving forces behind this creation was the quest for a ``Weil
cohomology theory'' of algebraic varieties in characteristic $p$
that could be used to prove Weil Conjectures. As a minimum
requirement, such a theory should have coefficients of
charactersitic $0$. An early attempt to construct it was due to
Serre \cite{S} who used exactly the functoriality of Witt
vectors. By functoriality, one can do the construction locally with
respect to the Zariski topology and endow an algebraic variety $X$
over a field of characteristic $p$ with the sheaf of rings
$W(\calo_X)$ that in good cases has no $p$-torsion. Taking its
cohomology groups and inverting $p$, one obtains a cohomology theory
with coefficients in characteristic $0$, and one can check whether
it has the other required properties.

Unfortunately, it does not. However, in an alternative universe, it
is pretty clear where the story would go. The reason Serre's
cohomology theory did not work was that it was not topological: it
was a characteristic-$0$ lifting of the cohomology groups
$H^\hdot(X,\calo_X)$ and not of say de Rham cohomology groups
$H^\hdot_{DR}(X)$. Thus a natural thing to do was to find a
functorial lifting $W\Omega_X^\hdot$ of the whole de Rham complex
$\Omega^\hdot_X$, or in other words, to add higher degree terms to
$W(O_X) = W\Omega^0_X$.

As we know, this is not what happened in real life. In the short
run, a much more successful suggestion for Weil cohomology appeared
in the form of \'etale cohomology of Grothendieck --- this has
coefficients in $\Q_l$, not $\Q_p$, but for Weil Conjectures, this
works just as well. Slightly later, Grothendieck turned his
attention to de Rham cohomology and developed his theory of
cristalline cohomology. For a variety in characteristic $p$, this
has coefficients in $\Z_p$, works for Weil Conjectures, and in fact
gives a functorial characteristic-$0$ lifting of de Rham cohomology
in a very precise sense. The price to pay is a very high-tech
definition of cristalline cohomology: it requires the full force of
topos theory, and does not provide cohomology groups in any explicit
form.

However, our story does not end here --- and for the next step, we
need to fast-forward twenty more years, to the late 70-ies. We stay
in France, and surprisingly, we rejoin the alternative universe:
using earlier work by Bloch \cite{Bl} and input from Deligne,
Illusie discovered that cristalline cohomology of a smooth algebraic
variety $X$ can in fact be computed by a functorial lifting
$W\Omega_X^\hdot$ of the de Rham complex $\Omega^\hdot_X$ to
characteristic $0$. This lifting is known as the {\em de Rham-Witt
  complex}.

What is remarkable here is that Bloch was in fact working with
algebraic $K$-theory, not cohomology, and $K$-theory exists in much
larger generality --- it is defined for any associative ring $A$,
with no commutativity assumptions. Moreover, it was known since 1962
that differential forms also exist in the same generality --- by a
theorem of Hochschild, Kostant and Rosenberg, differential forms of
degree $i$ on an affine smooth algebraic variety $X = \Spec A$ are
canonically identified with elements in the $i$-th Hochschild
homology group $HH_i(A)$, and Hochschild homology makes sense in the
non-commutative case. But in the 1970-ies, these ideas were pretty
far outside of mainstream --- in fact, Bloch was using Milnor
$K$-theory that requires commutativity, since good technology for
working with the full Quillen $K$-theory did not yet exist.

To get to the next step of the story, we have to skip twenty more
years that saw, among other things, the discovery of cyclic homology
by Connes and Tsygan, and the appearance of Non-Commutative Geometry
as a separate subject. We also have to move to Denmark, and change
the area from algebraic geometry to algebraic topology. A
topological version of Hochschild homology was introduced by
B\"okstedt in 1985, and in 1992, in a groundbreaking paper
\cite{BHM}, B\"okstedt, Hsiang and Madsen expanded it greatly in
their theory of Topological Cyclic Homology and cyclotomic
trace. Among a wealth of other results, the paper associates a
certain spectrum $TR(A,p)$ to any ring spectrum $A$ and prime
$p$. In particular, a ring $A$ is tautologically a ring spectrum, so
$TR(A,p)$ is well-defined. Moreover, if $A$ is annhilated by $p$,
then $TR(A,p)$ is an Eilenberg-Mac Lane spectrum, thus effectively a
complex of abelian groups. The homotopy groups $\pi_\idot TR(A,p)$
are the homology groups of this complex, and it has been shown by
Hesselholt \cite{hedRW}, \cite{heCo} that when $A$ is a commutative
algebra over a perfect ring $k$ of characteristic $p$ with smooth
spectrum $X = \Spec A$ --- that is, in the situation of the
Hochschild-Kostant-Rosenberg Thereom --- the groups $\pi_i TR(A,p)$
are naturally identified with the terms $H^0(X,W\Omega^i_X)$ of the
de Rham-Witt complex. In particular, in degree $0$, we recover the
$p$-typical Witt vectors $W(A)$.

For a general associative $A$, $\pi_0 TR(A,p)$ has been constructed
in a purely algebraic manner, again by Hesseholt \cite{hewi},
\cite{heerr}; he called it the {\em group of non-commutative Witt
  vectors $W(A)$}. For higher homotopy groups, a purely algebraic
construction is not known.

We are now running out of time to skip: fast-forwarding another
twenty years brings up more-or-less to today. What we want to
present in this paper, then, is a sort of a postscript to the long
story sketched above. For any associative unital algebra $A$ over a
perfect field $k$ of characteristic $p$, the {\em Hochschild-Witt
  homology groups} $WHH_\idot(A)$ have been constructed recently in
\cite{ka.witt}, \cite{ka.hw}. These ought to coincide with the
homotopy groups of the spectrum $TR(A,p)$ but at present, this has
not been checked. What has been proved is two comparison theorems:
one says that in the Hochschild-Kostant-Rosenberg situation,
$WHH_\idot(A)$ is identified with de Rham-Witt forms, and another
one says that in the general case, $WHH_0(A)$ is identified with
non-commutative Witt vectors $W(A)$ of Hesselholt.

\medskip

The paper is organized as follows. We start with a brief overview of
the classical story, presented in a certain way. This is
Section~\ref{cl.sec}. Section~\ref{nc.sec} shows that with very few
strategically placed modifications, the classical theory gives much
more. However, what it gives is not a theory for non-commutative
ring: rather, we obtain a notion of ``Witt vectors'' of a vector
space. Here we actually follow the approach sketched in
\cite{ka.trace} --- in order to generalize Hochschild homology, one
first has to treat it as a theory of two variables, an algebra $A$
over a field $k$ and an $A$-bimodule $M$, and work out what happens
when $A$ is just $k$ (but $M$ is arbirary). Then a general machine
of \cite{ka.trace} shows that as soon as the functor on vector
spaces is equipped with an additional structure, it automatically
extends to a functor on pairs $\langle A,M \rangle$. The structure
in question is that of a {\em trace functor}, and this is the
subject of Section~\ref{tr.sec}. In Section~\ref{hh.sec}, we turn to
homology: we construct our generalization of $HH_\idot(A)$ that we
call {\em Hochschild-Witt homology}, and we explain how it is
related to non-commutative Witt vectors of Hesselholt and to the de
Rham-Witt complex. Finally, in Section~\ref{add.sec}, we discuss
some parts of the theory that we skipped in Section~\ref{cl.sec}
(multiplication, Frobenius map), and discuss possible extensions of
the existing theory.

\section{Classical theory.}\label{cl.sec}

Fix a prime $p$. Denote by $R:\Z/p^{n+1}\Z \to \Z/p^n\Z$, $n \geq 1$,
the projection maps, and let
$$
\Z_p = \displaystyle\lim_{\overset{R}{\gets}}\Z/p^n\Z
$$
be the ring of $p$-adic integers. Alternatively, elements $a \in \Z_p$
can be thought of as formal power series
$$
a = \sum_{n \geq 0}a_ip^i,
$$
where $a_i$ lies in any fixed set $S \subset \Z$ of representatives
of residues $\overline{a} \in \F_p = \Z/p\Z$. For example, one can
take $S = \{0,1,\dots,p-1\}$, or $\{ 1,2,\dots,p\}$. There is
no canonical choice.

The key observation of Teichm\"uller is that we do have a canonical
choice if we allow $S$ to lie in $\Z_p$ itself, not in $\Z$. Namely,
one has the following.

\begin{lemma}\label{teich.le}
For any $x \in \Z/p\Z$, there exists a unique element $[x] \in \Z_p$
such that $[x]^p=[x]$ and $[x] = x \mod p$.
\end{lemma}

\proof{} Take $a \in \Z_p$ such that $a^p=a$. If $a$ is divisible by
$p$, then by induction, it is divisible by $p^n$ for any $n$, so
that $a=0$. If $a$ is not divisible by $p$, then for any $n \geq 0$,
it is invertible with respect to the multiplication in $\Z/p^n\Z$,
thus a root of unity of order $p-1$. All such roots form a subgroup
$\mu_{p-1}$ in the group $(\Z/p^n\Z)^*$ of invertible elements in
$\Z/p^n\Z$, and we need to prove that the projection
$R^{n-1}:\Z/p^n\Z \to \F_p$ induces an isomorphism $\mu_{p-1} \cong
\F_p^*$. This is clear, since $(\Z/p^n\Z)^*$ has order
$(p-1)p^{n-1}$, thus splits canonically into a product of $F_p^*$
and an abelian group of order $p^{n-1}$.
\endproof

The element $[x]$ of Lemma~\ref{teich.le} is known as the {\em
  Teichm\"uller representative} of the residue class $x$. By virtue
of uniqueness, the correspondence $x \mapsto [x]$ is multplicative,
thus defines a multiplicative map $T:\F_p \to \Z_p$ that splits the
projection $R:\Z_p \to \F_p$ (that is, we have $R \circ T =
\Id$). More generally, we obtain a functorial isomorphism
\begin{equation}\label{T.eq}
T_\idot:\prod_{i \geq 0}\F_p \cong \Z_p, \qquad T_\idot(\langle
a_0,a_1,\dots \rangle) = \sum_ip^i[a_i].
\end{equation}
This is only an isomorphism of sets. The question is, how do we
write down the ring operations in $\Z_p$ in terms of the components
$a_0,a_1,\dots$?

In fact, there are two questions here: how do we write down
addition, and then how do we write down multiplication. Since
the Teichm\"uller map is multiplicative, the second question is
easier, and we will concentrate on the first one.

To get a good answer to the question, we need to generalize
it. Namely, for any commutative ring $A$ and integer $n \geq 1$, denote by
$W_n(A)$ the $n$-fold product of copies of $A$ --- that is, the set
of collections $\langle a_0,\dots,a_{n-1}\rangle$ of elements
$a_\idot \in A$ (these are the ``Witt vectors''). What we want to
construct is a collection of abelian group structures on $W_n(A)$ that is
functorial with respect to $A$, and such that \eqref{T.eq} is a
group isomorphism when $A = \F_p$.

Moreover, this last condition can be also generalized so that it
makes sense for an arbitrary $A$. Namely, while the Teichm\"uller
map $T$ cannot be written down by an explicit formula, one can write
down explicitly the composition map
$$
\begin{CD}
\Z/p^{n+1}\Z @>{R^n}>> \F_p @>{T}>> \Z_p @>>> \Z/p^{n+1}\Z.
\end{CD}
$$
To do this, one proves the following effective version of
Lemma~\ref{teich.le}.

\begin{lemma}\label{teich.bis.le}
An element $x \in \Z/p^{n+1}\Z$ is a Teichm\"uller representative of
some residue class if and only if $x^{p^n}=x$. For any $x \in
\Z/p^{n+1}\Z$ with residue $\overline{x}=R^n(x) \in \F_p$, we have
$x^{p^n} = [\overline{x}] = T(R^n(x)) \mod p^{n+1}$.
\end{lemma}

\proof{} As in Lemma~\ref{teich.le}, use the canonical abelian group
decomposition $(\Z/p^n\Z)^* \cong \Z/p^n\Z \times \F_p^*$.
\endproof

Now for any commutative ring $A$ and any integer $n \geq 1$, define
the {\em ghost map} $w_n:W_{n+1}(A) \to A$ by
\begin{equation}\label{gh.eq}
w_n(\langle a_0,\dots,a_n \rangle) = \sum_{i=0}^np^ia_0^{p^{n-i}}.
\end{equation}
Since the restriction map $W_n(R^n):W_n(\Z/p^{n+1}\Z) \to W_n(\F_p)$
is surjective, the map $T_n:W_n(\F_p) \to \Z/p^{n+1}\Z$ induced by
\eqref{T.eq} is compatible with the group structure if and only if so
is the composition $T_n \circ W_n(R^n):W_n(\Z/p^{n+1}\Z) \to
\Z/p^{n+1}\Z$. By Lemma~\ref{teich.bis.le} and \eqref{T.eq}, this
composition is exactly the ghost map \eqref{gh.eq} for $A =
\Z/p^{n+1}\Z$. Thus if we require that our functorial group
structures on $W_\idot(A)$ make \eqref{gh.eq} additive for any $A$,
the map \eqref{T.eq} will automatically be an isomorphism of abelian
groups.

Let now $R:W_{n+1}(A) \to W_n(A)$, $V:W_n(A) \to W_{n+1}(A)$ be the
maps given by
$$
R(\langle a_0,\dots,a_n\rangle) = \langle a_0,\dots,a_{n-1} \rangle,
\quad V(\langle a_0,\dots,a_{n-1}\rangle) = \langle
0,a_0,\dots,a_{n-1} \rangle,
$$
where $V$ traditionally stands for {\em Verschiebung}, German for
``shift''. With this notation, Witt's theorem --- or rather, its
additive part --- can be formulated as follows.

\begin{theorem}\label{witt.thm}
There exists a unique set of functorial abelian group structures on
$W_n(A)$, $n \geq 1$, such that $R$, $V$ and $w_n$, $n \geq 1$ are
additive.
\end{theorem}

While this theorem can be proved directly, it is more instructive to
deduce it from the following elementary statement.

\begin{lemma}\label{c.le}
There exists a unique collection of polynomials $c_i(-,-)$, $i \geq
1$ with integer coefficients such that for any $n$ and commuting
$x_0$, $x_1$, we have
\begin{equation}\label{c.eq}
(x_0+x_1)^{p^n} = x_0^{p^n} + x_1^{p^n} +
\sum_{i=1}^np^ic_i(x_0,x_1)^{p^{n-i}}.
\end{equation}
\end{lemma}

\proof{} Uniqueness is obvious. Existence can be deduced from the
standard divisibility properties of binomial coefficients; we skip
this, since we will prove a more general statement later in
Proposition~\ref{c.prop}.
\endproof

\proof[Proof of Theorem~\ref{witt.thm}.] Since $w_0 = \id$, we must
have $W_1(A)=A$. Assume by induction that we have proved the
existence and uniqueness of the functorial group structures on
$W_i(A)$, $i \leq n$ with additive $R$, $V$ and $w_i$, $i \leq
n$. We then have to construct a functorial extension
$$
\begin{CD}
0 @>>> W_n(A) @>{V}>> W_{n+1}(A) @>{R^n}>> W_1(A)=A @>>> 0
\end{CD}
$$
of abelian groups, where $W_{n+1}(A) \cong A \times W_n(A)$ as
sets. With this decomposition, the operation in $W_{n+1}(A)$ must be
given by
\begin{equation}\label{c.w.eq}
\langle a_0,b_0 \rangle + \langle a_0,b_0 \rangle = \langle a_0 +
a_1,b_0 + b_1 + c_\idot(a_0,a_1) \rangle
\end{equation}
for some universal cocycle $c_\idot(-,-)$ that is symmetric with
respect to the transposition of variables. Moreover, since $c_\idot(-,-)$
is functorial in $A$, it is a map between affine schemes, thus it
must be polynomial in $a_0$ and $a_1$. Then $w_n$ is additive with
respect to the operation \eqref{c.w.eq} if and only if
$c_\idot(-,-)$ satisfies \eqref{c.eq}, so that Lemma~\ref{c.eq}
proves uniqueness. To prove existence, it suffices to define the
group operation by \eqref{c.w.eq} with the polynomials
$c_\idot(-,-)$ of Lemma~\ref{c.eq}, and check that it is indeed
associative and commutative. By functoriality, it suffices to check
this for $A=\Z$. But then the map $w_\idot = \langle w_0,\dots,w_n
\rangle:W_{n+1}(\Z) \to \Z^{n+1}$ is injective, so that the required
identities can be checked after composing with $w_\idot$, and they
follow from the additivity of $w_\idot$.
\endproof

We note that by construction, for any integers $m,n \geq 1$ and any
ring $A$, we have a functorial short exact sequence
\begin{equation}\label{R.V.seq}
\begin{CD}
0 @>>> W_m(A) @>{V^n}>> W_{m+n}(A) @>{R^m}>> W_n(A) @>>> 0
\end{CD}
\end{equation}
of abelian groups. In particular, we can take $m=1$; in this case,
we see that $W_{n+1}(A)$ is a functorial extension of $W_n(A)$ by
$A$.

\section{Non-commutative theory.}\label{nc.sec}

To adapt the construction of the Witt vectors $W(A)$ to a
non-commutative ring $A$, we start with Lemma~\ref{c.le}. Does the
recursive formula \eqref{c.le} admit a non-commutative refinement?
Naively, one would like to have non-commutative polynomials
$c_i(-,-)$, $i \geq 1$, of degrees $p^i$ such that for any $n$, we
have
\begin{equation}\label{c.nc.eq}
(x_0+x_1)^{\otimes p^n} = x_0^{\otimes p^n} + x_1^{\otimes p^n} +
\sum_{i=1}^np^ic_i(x_0,x_1)^{\otimes p^{n-i}} \in T^{p^n}(x_0,x_1),
\end{equation}
where $T^{p^n}(x_0,x_1)$ is the component of degree $p^n$ of the
free associative $\Z$-algebra $T^\hdot(x_0,x_1)$ on two variables
$x_0$, $x_1$. A moment's reflection shows that this is not possible
--- already for $n=1$, the non-commutative polynomial
$(x_0+x_1)^{\otimes p} - x_0^{\otimes p} - x_1^{\otimes p}$ is not
divisible by $p$. Something in \eqref{c.nc.eq} has to be
modified. This something turns out to be the coefficient $p^i$.

\begin{notation}\label{sigma.not}
For any free $\Z$-module $N$ and integer $i \geq 1$, we denote by
$\sigma:N^{\otimes i} \to N^{\otimes i}$ the permutation of order
$i$.
\end{notation}

In particular, for any $i$, we have $T^i(x_0,x_1) =
\Z[x_0,x_1]^{\otimes i}$, so that we have the permutation
$\sigma:T^i(x_0,x_1) \to T^i(x_0,x_1)$.

\begin{prop}\label{c.prop}
There exist a collection of non-commutative polynomials $c_i(-,-)$,
$i \geq 1$, of degree $p^i$ such that for any $n \geq 1$, we have
\begin{equation}\label{nc.eq}
(x_0+x_1)^{\otimes p^n} = x_0^{\otimes p^n} + x_1^{\otimes p^n} +
\sum_{i=1}^n(\id+\sigma+\dots+\sigma^{p^i-1})c_i(x_0,x_1)^{\otimes
  p^{n-i}}.
\end{equation}
\end{prop}

We note that Proposition~\ref{c.prop} is a refinement of the
existence part of Lemma~\ref{c.eq} (to obtain commutative
polynomials, just project to the symmetric algebra, and note that
$\sigma$ then intertwines the identity map). Just as
Lemma~\ref{c.le}, Proposition~\ref{c.prop} admits a direct
combinatorial proof, although it is rather non-trivial
(\cite[Proposition 1.2.3]{hewi} claims slightly less but actually
proves exactly \eqref{nc.eq}). Let us present an alternative
proof. It is somewhat roundabout but requires much less
computations, and gives useful by-products.

We start by observing that if for some integer $m \geq 1$, we have
already constructed the polynomials $c_i(-,-)$, $i < m$ such that
\eqref{nc.eq} holds for $n < m$, then to construct $c_m(-,-)$, we
need to shows that the difference
\begin{equation}\label{diff.eq}
(x_0+x_1)^{\otimes p^m} - x_0^{\otimes p^m} + x_1^{\otimes p^m} +
\sum_{i=1}^{m-1}(\id+\sigma+\dots+\sigma^{p^i-1})c_i(x_0,x_1)^{\otimes
  p^{m-i}}
\end{equation}
lies in the image of the map $\id + \sigma + \dots +
\sigma^{p^m-1}$. We note that this difference is obviously invariant
under $\sigma$.

\begin{defn}\label{Q.def}
For a free $\Z$-module $N$ and integer $n \geq 1$, $Q_n(N)$ is the
cokernel of the map
\begin{equation}\label{Q.eq}
\begin{CD}
\left(N^{\otimes p^n}\right)_\sigma @>{\id + \sigma + \dots +
  \sigma^{p^n-1}}>> \left(N^{\otimes p^n}\right)^\sigma.
\end{CD}
\end{equation}
\end{defn}

Equivalently, $Q_n(N) = \vH^0(\Z/p^n\Z,N^{\otimes p^n})$ is the
$0$-th Tate cohomology group of the cyclic group $\Z/p^n\Z$ acting
on $N^{\otimes p^n}$ via the permutation $\sigma$.

\begin{lemma}\label{T.le}
The correspondence $x \mapsto x^{\otimes p^n}$ factors as
$$
\begin{CD}
N @>>> N/pN @>{T_n}>> Q_n(N)
\end{CD}
$$
for some functorial map $T_n:N/pN \to Q_n(N)$.
\end{lemma}

\proof{} We need to show that for any $x_0,x_1 \in N$,
$(x_0+px_1)^{\otimes p^n} - x_0^{\otimes p^n}$ lies in the image of
the map $\id + \sigma + \dots + \sigma^{p^n-1}$. This is an
elementary direct computation (see \cite[Lemma 2.2]{ka.witt} for
details).
\endproof

\begin{corr}\label{W.corr}
There exists a unique functor $W_n$ from $\F_p$-vector spaces to
abelian groups such that for any free $\Z$-module $N$, we have
$$
Q_n(N) \cong W_n(N/pN),
$$
and this isomorphism is functorial in $N$.
\end{corr}

\proof{} Since every $\F_p$-vector space $M$ is of the form $M=N/pN$
for some free $\Z$-module $N$, defining $W_n$ on objects is trivial,
and the issue is morphisms: we need to show that $Q_n$ factors
through the full essentially surjective functor $N \mapsto
N/pN$. Explicitly, this means that for any two maps $a,a':N_0 \to
N_1$ between free $\Z$-modules, we have $Q_n(a) = Q_n(a + p
a')$. But by definition, the map $Q_n(a)$ is induced by $a^{\otimes
  p^n}:N_0^{\otimes p^n} \to N_1^{\otimes p^n}$, and similarly for
$Q_n(a+pa')$, so we are done by Lemma~\ref{T.le}.
\endproof

\begin{defn}\label{w.def}
For any $\F_p$-vector space $M$, the abelian group $W_n(M)$ is the
group of {\em $n$-truncated polynomial Witt vectors} of the vector
space $M$.
\end{defn}

The reason for our terminology is that polynomial Witt vectors of
Definition~\ref{w.def} behave similarly to the usual Witt
vectors. Indeed, Lemma~\ref{T.le} provides a canonical map $T_n:M
\to W_n(M)$, an analog of the Teichm\"uller representative
map. Moreover, for any free $\Z$-module $N$ and integer $n \geq 1$, we have
$N^{\otimes p^n} \cong (N^{\otimes p})^{\otimes p^{n-1}}$, and the
order-$p^n$ permutation $\sigma$ on $N^{\otimes p^n}$ induces an
order-$p$ automorphism $\overline{\sigma}$ on $Q_{n-1}(N^{\otimes
  p})$. This descends to a functorial automorphism
$\overline{\sigma}:W_{n-1}(M^{\otimes p}) \to W_{n-1}(M^{\otimes p})$
of order $p$, and we have a functorial map
\begin{equation}\label{V.eq}
V = \id + \overline{\sigma} + \dots +
\overline{\sigma}^{p-1}:W_{n-1}(M^{\otimes p}) \to W_n(M),
\end{equation}
an analog of the Verschiebung map. To complete the analogy, we need
to define restriction maps $R$. To do this, we need the following
version of Definition~\ref{Q.def}.

\begin{defn}\label{Q.bis.def}
For a free $\Z$-module $N$ and integer $n \geq 1$, $Q'_n(N)$ is the
cokernel of the map
$$
\begin{CD}
\left(N^{\otimes p^n}\right)_\sigma @>{\id + \sigma + \dots +
  \sigma^{p^{n+1}-1}}>> \left(N^{\otimes p^n}\right)^\sigma.
\end{CD}
$$
\end{defn}

\begin{lemma}\label{W.bis.le}
The correspondence $x \mapsto x^{p^n}$ factors through a functorial
map $T'_n:N/pN \to Q'_n(N)$. Moreover, there exists a unique
functor $W'_n$ from $\F_p$-vector spaces to abelian group such that
$Q'_n(N) \cong W'_n(N/pN)$ functorially in $N$, and $T'_n$ descends
to a functorial map $T'_n:M \to W'_n(M)$ for any $\F_p$-vector
space $M$.
\end{lemma}

\proof{} Same as Lemma~\ref{T.le} and Corollary~\ref{W.corr} (for
details, see \cite[Lemma 2.2]{ka.witt}).
\endproof

On the other hand, the permutation $\sigma:N^{\otimes p^n} \to
N^{\otimes p^n}$ has order $p^n$, so that
$$
\id + \sigma + \dots + \sigma^{p^{n+1}} = p(\id + \sigma + \dots +
\sigma^{p^n}),
$$
and we have a natural functorial projection $W'_n(M) \to W_n(M)$.

\begin{defn}
For any free $\Z$-module $N$, an additive map $c:N \to
\left(N^{\otimes p}\right)^\sigma$ is {\em standard} if the induced
map $\overline{c}:N/pN \to Q_1(N)$ is the Teichm\"uller map $T_1$ of
Lemma~\ref{T.le}.
\end{defn}

Standard maps exist (e.g.\ choose a basis in $N$, and let $c$ send
any basis element $e$ to $e^{\otimes p}$) and enjoy the following
property.

\begin{prop}\label{stand.prop}
For any standard map $c:N \to N^{\otimes p}$ and any integer $n \geq
1$, the induced map $c^{\otimes p^n}:Q'_n(N) \to Q_{n+1}(N)$ is an
isomorphism that does not depend on the choice of $c$, and we have
$c^{\otimes p^n} \circ T'_n = T_{n+1}$.
\end{prop}

\proof{} See \cite[Proposition 2.6]{ka.witt}. \endproof

As a corollary of Proposition~\ref{stand.prop}, we have a functorial
isomorphism $W'_n(M) \cong W_{n+1}(M)$ for any $\F_p$-vector space
$M$. Composing the inverse isomorphism with the projection $W'_n(M)
\to W_n(M)$, we obtain a functorial map
$$
R:W_{n+1}(M) \to W_n(M),
$$
an analog of the restriction map. Interaction between the maps $R$
and $V$ is also analogous to what one has for classical Witt vectors
--- in particular, it has been proved in \cite[Lemma 3.7]{ka.witt}
that for any $n,m \geq 0$ and $\F_p$-vector space $M$, we have an
exact sequence
\begin{equation}\label{R.V.nc.seq}
\begin{CD}
W_m(M^{\otimes p^n}) @>{V^n}>> W_{m+n}(M) @>{R^m}>> W_n(M) @>>> 0
\end{CD}
\end{equation}
of abelian groups, a vector space counterpart of the sequence
\eqref{R.V.seq}.

\proof[Proof of Proposition~\ref{c.prop}.] Assume by induction that
we have defined the polynomials $c_i(-,-)$ for $i \leq n-1$. Then
projecting the equality \eqref{nc.eq} in $T^{p^{n-1}}(x_0,x_1)$ to
$Q'_{n-1}(\Z[x_0,x_1])$, we see that
$$
T'_{n-1}(x_0+x_1) = T'_{n-1}(x_0) + T'_{n-1}(x_1) +
\sum_{i=1}^{n-1}V^i(T'_{n-1-i}(c_i(x_0,x_1))).
$$
However, $Q'_{n-1}(\Z[x_0,x_1]) \cong Q_n(\Z[x_0,x_1])$ by
Proposition~\ref{stand.prop}, and in terms of the group
$Q_n(\Z[x_0,x_1])$, this equality reads as
$$
(x_0+x_1)^{\otimes p^n} = x_0^{\otimes p^n} + x_1^{\otimes p^n} +
\sum_{i=1}^{n-1}(\id+\sigma+\dots+\sigma^{p^i-1})c_i(x_0,x_1)^{\otimes
  p^{n-i}}.
$$
In other words, \eqref{diff.eq} vanishes after projecting to
$Q_n(\Z[x_0,x_1])$. As we have remarked, this implies the existence
of $c_n(-,-)$ and proves the claim.
\endproof

\section{Trace functors.}\label{tr.sec}

\subsection{Trace maps.}

One obvious difference between Lemma~\ref{c.le} and
Proposition~\ref{c.prop} is that the latter does not claim
uniqueness. Uniqueness is in fact wrong; formally, this is reflected
in the fact that the sequence \eqref{R.V.nc.seq} is not exact on the
left. If $m=1$, then the image of the map $V^n:W_1(M^{\otimes p^n})
\to W_{n+1}(M)$ is the space $M^{\otimes p^n}_\sigma$ of
coinvariants with respect to the permutation $\sigma$. At the
opposite extreme, if $n=1$, then the kernel of the map
$R^n:W_{n+1}(M) \to W_1(M)=M$ is the group of coinvariants
$W_n(M^{\otimes p})_{\overline{\sigma}}$ with respect to the
permutation $\overline{\sigma}$ of \eqref{V.eq}. To understand the
behaviour in general, we observe that $W_\idot(-)$ carry an
additional structure: that of a {\em trace functor} in the sense of
\cite[Section 2]{ka.trace}

\begin{defn}\label{tr.def}
A {\em trace functor} from a unital monoidal category $\C$ to a
category $\E$ is a collection of a functor $F:\C \to \E$ and
isomorphisms
\begin{equation}\label{M.N}
\tau_{M,N}:F(M \otimes N) \cong F(N \otimes M)
\end{equation}
for any $M,N \in \C$ such that $\tau_{M,N}$ are functorial in $M$
and $N$, $\tau_{1,M} = \id$ for any object $M \in \C$, and for any
three objects $M,N,L \in \C$, we have
\begin{equation}\label{M.N.L}
\tau_{L,M,N} \circ \tau_{N,L,M} \circ \tau_{M,N,L} = \id,
\end{equation}
where $\tau_{A,B,C}$ for any $A,B,C \in \C$ is the composition of
the map $\tau_{A,B \otimes C}$ and the associativity isomorphism $(B
\otimes C) \otimes A \cong B \otimes (C \otimes A)$.
\end{defn}

If a monoidal category $\C$ is symmetric, such as for example the
category of vector spaces over a field, then any functor $\C \to \E$
is tautologically a trace functor, with $\tau_{-,-}$ induced by the
commutativity isomorphism in $\C$. However, even in this case, there
are also non-trivial trace functor structures. The basic example
considered in \cite{ka.trace} is the following.

\begin{exa}
Let $\C$ be the category of free $\Z$-modules, and fix an integer $l
\geq 2$. For any $M \in \C$, let $\sigma_M=\sigma:M^{\otimes l} \to
M^{\otimes l}$ be as in Notation~\ref{sigma.not}, and for any $M,N
\in \C$, let
\begin{equation}\label{tau.eq}
\wt{\tau}_{M,N}:(M \otimes N)^{\otimes l} \to (N \otimes M)^{\otimes
  l}
\end{equation}
be the product of the $l$-th tensor power of the commutativity
morphism and the automorphism $\sigma_N \otimes \id_{M^{\otimes l}}$
of $(N \otimes M)^{\otimes l} \cong N^{\otimes l} \otimes M^{\otimes
  l}$. Then $\wt{\tau}_{M,N}$ obviously commutes with $\sigma_{N
  \otimes M}$, thus induces a map
$$
\tau_{M,N}:(M \otimes N)^{\otimes l}_\sigma \to (N \otimes
M)^{\otimes l}_\sigma
$$
on the modules of coinvariants with respect to $\sigma$. Setting
$C_l(M) = M^{\otimes l}_\sigma$ with these structure maps
$\tau_{M,N}$ gives a trace functor from $\C$ to itself.
\end{exa}

Alternatively, one can consider the module $C^l(M) =
\left(M^{\otimes l}\right)^\sigma$ of invariants with respect to
$\sigma$; then \eqref{tau.eq} turns $C^l(M)$ into a trace functor as
well. Moreover, if $l=p^n$, then the map \eqref{Q.eq} is a map of
trace functors, so that the functor $Q_n$ of
Definition~\ref{Q.def} is naturally a trace functor from free
$\Z$-modules to abelian groups. This trace functor structure then
descends to a trace functor structure on the Witt vectors functor
$W_n(-)$.

It is not difficult to deduce from \eqref{M.N.L} that for any trace
functor $F:\C \to \E$ and for any $l$ objects $M_1,\dots,M_l \in
\C$, we have a natural map
\begin{equation}\label{tau.M}
\tau:F(M_1 \otimes \dots \otimes M_l) \to F(M_2 \otimes \dots
\otimes M_l \otimes M_1),
\end{equation}
and these maps satify an $l$-variable version of \eqref{M.N.L}. In
particular, for any $M \in \C$, we have an automorphism
$\tau:F(M^{\otimes l}) \to F(M^{\otimes l})$ of order $l$. Then it
has been proved in \cite[Lemma 3.7]{ka.witt} that for any $n,m \geq
1$, the sequence \eqref{R.V.nc.seq} can be refined to a functorial
short exact sequence
\begin{equation}\label{W.seq}
\begin{CD}
0 @>>> W_m(M^{\otimes p^n})_\tau @>>> W_{m+n}(M) @>>> W_n(M) @>>> 0,
\end{CD}
\end{equation}
where $\tau:W_m(M^{\otimes p^n}) \to W_m(M^{\otimes p^n})$ is the
automorphism of order $p^n$ induced by the trace functor structure
on $W_m$. If $m=1$, then $\tau=\sigma$ is the permutation, and if
$n=1$, then $\tau=\overline{\sigma}$ is as in \eqref{V.eq}. We
caution the reader that for $m \geq 2$, the trace functor structure
on $W_m$ is non-trivial, so that $\tau$ is different from the map
induced by the permutation $\sigma:M^{\otimes p^n} \to M^{\otimes
  p^n}$.

\subsection{Small categories.}

While constructing the map \eqref{tau.M} directly is not difficult,
this has not been done in \cite{ka.trace}. Instead, the maps are
deduced from a convenient packaging of a trace functor structure
using the category $\Lambda$ introduced by Connes in
\cite{connes}. This is a small category whose objects $[n]$ are
indexed by positive integers $n$. Maps between $[n]$ and $[m]$ can
be defined in various equivalent ways, see e.g.\ \cite{Lo} or
\cite{FTadd}; for the convenience of the reader, we recall the most
geometric of the descriptions.
\begin{itemize}
\item The object $[n] \in \Lambda$ is thought of as a ``wheel'' -- a
  cellular decomposition of the unit circle $S^1 \subset \C$ with
  $n$ vertices corresponding to $n$-th roots of unity, and $n$
  edges. A continuous map $f:[n'] \to [n]$ is {\em good} if it is
  cellular, has degree $1$, and the induced map $\wt{f}:\R \to \R$
  between the universal covers is non-decreasing. Morphisms from
  $[n']$ to $[n]$ in the category $\Lambda$ are homotopy
  classes of good maps $f:[n'] \to [n]$.
\end{itemize}
We will denote the set of vertices of $[n] \in \Lambda$ by
$V([n])$. Every map $f:[n'] \to [n]$ in $\Lambda$ induces a map
$V(f):V([n']) \to V([n])$, and for any $v \in V([n])$, the preimage
$V(f)^{-1}(v) \subset V([n'])$ has a natural total order induced by
the clockwise orientation of the circle.

Now, for any unital monoidal category $\C$, define a category
$\C^\hush$ as follows:
\begin{itemize}
\item objects are pairs $\langle [n],\{c_\idot\} \rangle$ of an
  object $[n] \in \Lambda$ and a collection of objects $c_v \in \C$
  indexed by vertices $v \in V([n])$,
\item morphisms from $\langle [n],c_\idot \rangle$ to $\langle
  [n'],c'_\idot \rangle$ are given by pairs $\langle f,\{ f_v \}
  \rangle$ of a morphism $f:[n] \to [n']$ and a collection of
  morphisms
$$
f_v:\otimes_{v' \in V(f)^{-1}(v)}c_{v'} \to c_v, \qquad v \in
V([n]),
$$
where the product is taken in the natural order on $V(f)^{-1}(v)$.
\end{itemize}
A morphism $\langle f,\{f_v\} \rangle$ in $\C^\hush$ is called {\em
  cartesian} if all the maps $f_v$ are invertible. With these
definitions, one has the following result.

\begin{lemma}[{{\cite[Lemma 2.3]{ka.trace}}}]\label{tr.le}
Giving a trace functor from $\C$ to some category $\E$ is equivalent
to giving a functor $\C^\hush \to \E$ that sends any cartesian map
to an invertible map.\endproof
\end{lemma}

Here is a sketch of the correspondence of Lemma~\ref{tr.le}. The
category $\C$ is naturally embedded into $\C_\hush$ by sending $c
\in \C$ to $\langle [1],c \rangle$. Thus any functor
$F^\hush:\C^\hush \to \E$ gives by restriction a functor $F:\C \to
\E$. To see the map \eqref{M.N}, consider the object $\langle
  [2],\{M,N\}\rangle \in \C^\hush$, and note that there are two maps
  $s,t:[2] \to [1]$ in $\Lambda$ that impose opposite orders on the
  set $V([2])$. These maps lift to cartesian maps
$$
\wt{s}:\langle [2],\{M,N\} \rangle \to \langle [1],M \otimes N
\rangle, \qquad \wt{t}:\langle [2],\{M,N\} \rangle \to \langle [1],N
\otimes M\rangle,
$$
and we then have $\tau_{M,N} = F^\hush(\wt{t}) \circ
F^\hush(\wt{s})^{-1}$. The constraint \eqref{M.N.L} is encoded in
the structure of the category $\C^\hush$, and so is the map
\eqref{tau.M} --- in effect, the group $\Aut([n])$ of automorphisms
of the object $[n] \in \Lambda$ is the cyclic group $\Z/n\Z$ whose
generator $\tau$ is the clockwise rotation by $2\pi/n$, and it lifts
to a map $\wt{\tau}$ in $\C^\hush$ that gives \eqref{tau.M} after
applying $F^\hush$.

\subsection{Fibrations.}

To describe the functors $W_n^\hush$ corresponding to the trace
functor structures on polynomial Witt vectors $W_n$, we need to make
a digression on homology of small categories. Recall that for any
small category $I$ and ring $k$, the category $\Fun(I,k)$ of
functors from $I$ to $k$-modules is abelian. For any functor
$\gamma:I \to I'$ between small categories, we have the natural
pullback functor $\gamma^*:\Fun(I',k) \to \Fun(I,k)$, $E \mapsto E
\circ \gamma$, and it has a left and right-adjoint {\em Kan
  extension functor} $\gamma_!,\gamma_*:\Fun(I,k) \to
\Fun(I',k)$. For example, if $I' = \ppt$ is the point category, then
$\gamma_! = \colim_I$ and $\gamma_*=\lim_I$ are the colimit and
limit functor; furthermore, if $I = \ppt_G$ is the groupoid with one
object with automorphism group $G$, then $\Fun(I,k)$ is the category
of representations of $G$ in $k$-modules, and $\gamma_!(E) = E_G$,
$\gamma^*(E) = E^G$ send a representation $E$ to its module of
$G$-coinvariants resp.\ $G$-invariants. For a more general target
category $I'$, computing Kan extensions can be cumbersome, but there
is one situation where it is still easy.

\begin{defn}\label{fib.def}
A functor $\pi:I_0 \to I_1$ is a {\em fibration in groupoids} if
\begin{enumerate}
\item for any morphism $f_1:i'_1 \to i_1$ in $I_1$ and any object $i_0
  \in I_0$ with $\pi(i_0)=i_1$, there exists a morphism $f_0:i'_0 \to
  i_0$ in $I_0$ such that $\pi(f_0)=f_1$, and
\item for any two such morphisms $f'_0:i'_0 \to i_0$, $f''_0:i''_0
  \to i_0$, there exsits a unique morphism $g:i'_0 \to i''_0$ such
  that $\pi(g)=\id_{i'}$ and $f'_0 = f''_0 \circ g$.
\end{enumerate}
A functor $\pi$ is a {\em bifibration in groupoids} if both $\pi$
and the opposite functor $\pi^o:I_0^o \to I_1^o$ are fibrations in
groupoids.
\end{defn}

Definition~\ref{fib.def} is a special case of the general formalism
of \cite{sga} but we will not need the full generality --- it
suffices to know that for any bifibration $\pi:I' \to I$ and any
morphism $f:i' \to i$, there exist a pair of adjoint functors
$f^*:I'_i \to I'_{i'}$, $f_!:I'_{i'} \to I'_i$ between the fibers
$I'_i=\pi^{-1}(i)$, $I'_{i'}=\pi^{-1}(i')$ of the bifibration
$\pi$. If $\pi$ is a bifibration in groupoids in the sense of
Definition~\ref{fib.def}, one easily checks that the fibers
$I'_{i'}$, $I'_i$ are groupoids, so that $f^*$ and $f_!$ are
mutually inverse equivalences of categories. We also need the
following special case of \cite[Lemma 1.7]{ka0}.

\begin{lemma}\label{bc.le}
Assume given a bifibration in groupoids $\pi:I_1 \to I_0$ with small
$I_0$, $I_1$, and a functor $\gamma_0:I'_0 \to I_0$ from some small
category $I'_0$. Consider the pullback square
$$
\begin{CD}
I_1' @>{\gamma_1}>> I_1\\
@V{\pi'}VV @VV{\pi}V\\
I_0' @>{\gamma_0}>> I_0.
\end{CD}
$$
Then $\pi'$ is a bifibration in groupoids, and the natural maps
$$
\pi'_! \circ \gamma_1^* \to \gamma_0^* \circ \pi_!, \qquad
\gamma_0^* \circ \pi_* \to \pi'_* \circ \gamma_1^*
$$
adjoint to the isomorphism $\gamma_1^* \circ \pi^* \cong \pi^{'*}
\circ \gamma_0^*$ are themselves invertible.\endproof
\end{lemma}

As a corollary of Lemma~\ref{bc.le}, we see that if we are given a
category $I'$, a functor $E:I' \to k\amod$ to the category of
modules over a ring $k$, and a bifibration in groupoids $\pi:I' \to
I$ with small fibers, then $\pi_!E$ and $\pi_*E$ are well-defined
even if $I'$ and $I$ are large --- indeed, for any small category
$I_0 \subset I$ with preimage $I'_0 = \pi^{-1}(I_0) \subset I'$, the
induced functor $\pi:I'_0 \to I_0$ is a bifibration in groupoids, so
that $\pi_!E$ and $\pi_*E$ are well-defined on $I_0$, and both do
not depend on the choice of $I_0 \subset I$. Moreover, applying
Lemma~\ref{bc.le} to the embedding $\gamma:\ppt \to I$ onto an
object $i \in I$, we obtain canonical identifications
$$
\pi_!E(i) \cong \colim_{I'_i}E, \qquad \pi_*E(i) \cong
\textstyle\lim_{I'_i}E,
$$
where as before, $I'_i = \pi^{-1}(i) \subset I'$ is the fiber of the
bifibration $\pi$.

\subsection{Edgewise subdivision.}

A typical example of the situation of Definition~\ref{fib.def}
occurs in the study of Connes' category $\Lambda$. Namely, fix an
integer $l \geq 1$, and note that for any $n \geq 1$, the object
$[nl] \in \Lambda$ has an automorphism $\sigma = \tau^n:[nl] \to
[nl]$ of order $l$ given by clockwise rotation by $2\pi/l$. Let
$\Lambda_l$ be the category whose objects $[n]$ correspond to
integers $n \geq 1$, and with maps from $[n]$ to $[n']$ given by
$\sigma$-equivariant maps from $[nl]$ to $[n'l]$ in $\Lambda$. We
then have a tautological embedding $i_l:\Lambda_l \to \Lambda$, $[n]
\mapsto [nl]$, known as the {\em edgewise subdivision functor}. On
the other hand, we can identify the quotient $S^1/\sigma$ with $S^1$
by the $l$-th power map $z \mapsto z^l$, and this identification is
compatible with our cellular decompositions. Then sending a map in
$\Lambda_l$ to the induced map of the quotient circles defines a
functor $\pi_l:\Lambda_l \to \Lambda$. This functor is a bifibration
in groupoids in the sense of Definition~\ref{fib.def}. Its fibers
are naturally identified with the groupoid $\ppt_l = \ppt_{\Z/l\Z}$
with one object with automorphism group $\Z/l\Z$.

More generally, assume given a unital monoidal category $\C$,
consider the category $\C^\hush$ of Lemma~\ref{tr.le}, and define a
category $\C^\hush_l$ by the fibered product square
$$
\begin{CD}
\C^\hush_l @>{\pi_l}>> \C^\hush\\
@VVV @VVV\\
\Lambda_l @>{\pi_l}>> \Lambda.
\end{CD}
$$
Then $\pi_l:\C^\hush_l \to \C^\hush$ is also a bifibration in
groupoids with fiber $\ppt_l$. Moreover, the edgewise subdivision
functor $i_l:\Lambda_l \to \Lambda$ extends to a functor
$i_l:\C^\hush_l \to \C^\hush$ sending $\langle [n],\{c_\idot\}
\rangle$ to $[nl]=i_l([n])$ with the collection of objects
$\{c_{q(v)}\}$, where $q:V([nl]) \to V([n])$ is induced by the
quotient map $q:S^1 \to S^1/\sigma$.

Now, if we take as $\C$ the category of free $\Z$-modules, then the
tautological embedding from $\C$ to the category $\Ab$ of abelian
groups has the trivial structure of a trace functor, thus defines a
functor $I^\hush:\C^\hush \to \Ab$. Then Lemma~\ref{bc.le} shows
that for any $l \geq 1$, both functors
$\pi_{l!}i_l^*I^\hush,\pi_{l*}i_l^*I^\hush:\C^\hush \to \Ab$ send
cartesian maps to invertible maps, thus correspond to trace functors
by Lemma~\ref{tr.le}. Moreover, the map $\tr_l = \id + \sigma +
\dots + \sigma^l$ defines a map of functors
\begin{equation}\label{tr.l}
\tr_l:\pi_{l!}i_l^*I^\hush \to \pi_{l*}i_l^*I^\hush,
\end{equation}
and if we take $l=p^n$, then the cokernel $Q_n^\hush$ of the map
\eqref{tr.l} corresponds to the trace functor structure on the
functor $Q_n$ of Definition~\ref{Q.def}. By virtue of
Corollary~\ref{W.corr}, $Q_n^\hush$ actually factors through a
functor $W_n^\hush:\overline{\C}^\hush \to \Ab$, where
$\overline{\C}$ is the category of $\F_p$-vector spaces; the functor
$W_n^\hush$ then corresponds to the trace functor structure on the
Witt vectors functor $W_n$. One then shows (see \cite[Proposition
  4.3]{ka.witt}) that the restriction maps $R$ extend to maps
$R:W_{n+1}^\hush \to W_n^\hush$, and the Verschiebung maps $V$
induce maps
\begin{equation}\label{V.la}
\pi_{p!}i_p^*W_n^\hush \to W_{n+1}^\hush.
\end{equation}
More generally, for any $n,m \geq 0$, we have a short exact sequence
\begin{equation}\label{W.tr.seq}
\begin{CD}
0 @>>> \pi_{p^n!}i_{p^n}^*W_n^\hush @>{V^m}>> W_{n+m}^\hush @>{R^n}>>
W_m^\hush @>>> 0
\end{CD}
\end{equation}
of functors from $\overline{\C}^\hush$ to $\Ab$, a further
refinement of the sequence \eqref{W.seq}.

\section{Hochschild homology.}\label{hh.sec}

\subsection{Generalities on cyclic homology.}

While the category $\Lambda$ does package nicely trace functors of
Definition~\ref{tr.def}, this was not why it was originally
introduced in \cite{connes} --- the intended application was to
cyclic homology.

To understand this, recall that for any small category $I$ and ring
$k$, the category $\Fun(I,k)$ of functors from $I$ to $k\amod$ is
abelian, with enough injectives and projectives. For any $E \in
\Fun(I,k)$, the homology and cohomology modules of $I$ with
coefficients in $E$ are given by
$$
H_\idot(I,E) = L^\hdot\colim_IE, \qquad H^\hdot(I,E) =
R^\hdot\textstyle\lim_IE.
$$
If $E = k$ is the constant functor with value $k$, and $k$ is
commutative, then $H^\hdot(I,k)$ is an algebra, and for any $E \in
\Fun(I,k)$, both $H_\idot(I,E)$ and $H^\hdot(I,E)$ are modules over
$H^\hdot(I,k)$.

\begin{exa}
Let $\Delta$ be the category of finite non-empty totally ordered
sets $[n]=\{0,\dots,n\}$ and order-preserving maps, and let
$\Delta^o$ be the opposite category. Then $\Fun(\Delta^o,k)$ is the
category of simplicial $k$-modules, and for any $E \in
\Fun(\Delta^o,k)$, the homology modules $H_\idot(\Delta^o,E)$ can be
computed by the standard chain complex of the simplicial $k$-module.
\end{exa}

Consider now the category $\Lambda$, and let $\Lambda/[1]$ be the
category of objects $[n] \in \Lambda$ equipped with a morphism
$f:[n] \to [1]$. Then $V([1])$ consists of a single vertex $o$, and
the preimage $V(f)^{-1}(o)$ carries a natural total order, so that we
can define a functor $\Lambda/[1] \to \Delta$ by sending $f:[n] \to
[1]$ to $V(f)^{-1}(o)$. It turns out that this is an equivalence of
categories. Composing the inverse equivalence with the forgetful
functor $\Lambda/[1] \to \Lambda$ sending $f:[n] \to [1]$ to $[n]$,
we obtain a natural embedding $j:\Delta \to \Lambda$. Moreover,
$\Lambda$ is self-dual: we have an equivalence $\Lambda \cong
\Lambda^o$ sending a cellular decomposition of $S^1$ to the dual
cellular decomposition. Therefore we also have an embedding
$j^o:\Delta^o \to \Lambda$.

\begin{defn}
For any ring $k$ and functor $E \in \Fun(\Lambda,k)$, the {\em
  Hochschild} and {\em cyclic homology groups} of $E$ are given by
$$
HH_\idot(E) = H_\idot(\Delta^o,j^{o*}E), \qquad
HC_\idot(E)=H_\idot(\Lambda,E).
$$
\end{defn}

It turns out (see e.g.\ \cite[Chapter 6]{Lo}) that the cohomology
$H^\hdot(\Lambda,k)$ with coefficients in a commutative ring $k$ is
naturally identified with the polynomial algebra $k[u]$ in one
generator $u$ of degree $2$, so that for any $E \in
\Fun(\Lambda,k)$, the cyclic homology $HC_\idot(E)$ is a
$k[u]$-algebra. One has a $k[u]$-equivariant spectral sequence
\begin{equation}\label{hdr.seq}
HH_\idot(E)[u^{-1}] \Rightarrow HC_\idot(E),
\end{equation}
known as {\em Hodge-to-de Rham} spectral sequence. One can also
invert the generator $u$ and define {\em periodic cyclic homology}
$HP_\idot(E)$ by
$$
HP_\idot(E) = R^\hdot\lim_{\overset{u}{\gets}}HC_\idot(E),
$$
where $R^\hdot\lim$ is the derived functor of the inverse limit
functor. Then the spectral sequence \eqref{hdr.seq} gives rise to a
spectral sequence
\begin{equation}\label{hdr.bis.seq}
HH_\idot(E)((u)) \Rightarrow HP_\idot(E),
\end{equation}
where on the left-hand side, we have formal Laurent power series in
the generator $u$.

Now assume given an associative unital $k$-algebra $A$, and consider
the functor $A^\hush \in \Fun(\Lambda,k)$ that sends an object $[n]
\in \Lambda$ to $A^{\otimes n}$, with copies numbered by vertices $v
\in V([n])$, and sends a morphism $f:[n'] \to [n]$ to the map
$$
A^\hush(f) = \bigotimes_{v \in V([n])} m_{V(f)^{-1}(v)},
$$
where for any finite totally ordered set $S$, $m_S:A^{\otimes S} \to A$ is
the product map. The Hochschild and cyclic homology of the algebra
$A$ are then given by
\begin{equation}\label{hh.hc.eq}
HH_\idot(A)=HH_\idot(A^\hush), \ HC_\idot(A)=HC_\idot(A^\hush),
\ HP_\idot(A)=HP_\idot(A^\hush),
\end{equation}
and one has the following comparison result (\cite[Chapter 3]{Lo},
\cite{FTadd}).

\begin{theorem}\label{hkr.thm}
Assume that $k$-algebra $A$ is flat, commutative and finitely
generated, and $X = \Spec A$ is smooth over $k$. Then we have the
canonical identifications
\begin{equation}\label{hkr.eq}
HH_i(A) \cong H^0(X,\Omega^i_X)
\end{equation}
for any $i \geq 0$, and under these identifications, the first
differential $B:HH_i(A) \to HH_{i+1}(A)$ in the spectral sequences
\eqref{hdr.seq}, \eqref{hdr.bis.seq} is identified with the de Rham
differential $d$.\endproof
\end{theorem}

The identification \eqref{hkr.eq} is the famous theorem of
Hochschild, Kostant and Rosenberg \cite{hkr} (they used the original
definition of $HH_\idot(A)$ obtained by computing
$HH_\idot(j^{o*}A^\hush)$ in \eqref{hh.hc.eq} via the standard
complex of the simplicial $k$-module $j^{o*}A^\hush$). The
differential $B$ was discovered by G. Rinehart \cite{rine},
forgotten and then rediscovered as a part of the cyclic homology
package, independently by A. Connes \cite{connes1} and B. Tsygan
\cite{tsy}, in about 1982. The packaging using the category
$\Lambda$ has been discovered by Connes one year later. The name
``Hodge-to-de Rham spectral sequence'' for the sequences
\eqref{hdr.seq}, \eqref{hdr.bis.seq} is motivated by
Theorem~\ref{hkr.thm}. If $k$ contains $\Q$, one can prove more: the
Hodge-to-de Rham spectral sequences degenerate at second term, and in
particular, one has $HP_\idot(A) \cong H^\hdot_{DR}(X)((u))$, the de
Rham cohomology of the variety $X = \Spec A$. The Theorem itself
works over any $k$, but even if $k$ is a field not containing $\Q$,
it is currently unknown whether \eqref{hdr.bis.seq} degenerates
(unless $A$ is a polynomial algebra, see \cite[Chapter 3.2]{Lo}).

\subsection{Twisting.}

We now observe that Lemma~\ref{tr.le} allows us to generalize
\eqref{hh.hc.eq} in the following way. Assume given a unital
monoidal category $\C$, and consider the category $\C^\hush$ with
the forgetful functor $\rho:\C^\hush \to \Lambda$. Then an
associative unital algebra object $A$ in $\C$ gives rise to a
section $\alpha:\Lambda \to \C^\hush$ of the projection $\rho$
sending $[n]$ to $\langle [n],\{A\} \rangle$, the collection of
copies of the object $A$ numbered by vertices $v \in V([n])$. If we
also have a trace functor $F$ from $\C$ to $k\amod$ for some ring
$k$, then we can consider the corresponding functor
$F^\hush:\C^\hush \to k\amod$. Composing $\alpha$ and $F^\hush$, we
obtain a functor $FA^\hush = F^\hush \circ \alpha \in
\Fun(\Lambda,k)$.

\begin{defn}
The {\em twisted Hochschild and periodic cyclic homology} of the
algebra object $A$ with respect to the trace functor $F$ is given by
\begin{equation}\label{hh.hc.tw.eq}
FHH_\idot(A)=HH_\idot(FA^\hush), \qquad FHP_\idot(A)=HP_\idot(FA^\hush).
\end{equation}
\end{defn}

In particular, we can take $k=\Z$, let $\C$ be the category of
$\F_p$-vector spaces, and take the polynomial Witt vectors functor
$W_n$ for some integer $n$. Then for any associative unital
$\F_p$-algebra $A$, we obtain the {\em Hochschild-Witt homology
  groups} $W_nHH_\idot(A)$ and corresponding periodic cyclic groups
$W_nHP_\idot(A)$. The restriction maps $R$ induce functorial maps
\begin{equation}\label{R.hw}
R:W_{n+1}HH_\idot(A) \to W_nHH_\idot(A),
\end{equation}
so that the groups $W_nHH_\idot(A)$, $n \geq 1$, form a projective
system. For any $n,m \geq 0$, the exact sequence \eqref{W.tr.seq}
induces an exact sequence
\begin{equation}\label{W.la.seq}
\begin{CD}
0 @>>> \pi_{p^n!}i_{p^n}^*W_nA^\hush @>{V^m}>> W_{n+m}W^\hush
@>{R^n}>> W_mA^\hush @>>> 0
\end{CD}
\end{equation}
in $\Fun(\Lambda,\Z)$. Moreover, for any integer $l \geq 1$, the
embedding $j^o:\Delta^o \to \Lambda$ lifts canonically to an
embedding $j^o_l:\Delta \to \Lambda_l$ such that $\pi_l \circ j^o_l
\cong j^o$, and $i_l \circ j^o_l = j^o \circ \overline{i}_l$ for a
certain functor $\overline{i}_l:\Delta^o \to \Delta^o$. The classic
Edgewise Subdivision Lemma \cite{edge} then claims, among other
things, that for any ring $k$ and object $E \in \Fun(\Delta^o,k)$,
the natural map
$$
H_\idot(\Delta^o,\overline{i}_l^*E) \to H_\idot(\Delta^o,E)
$$
is an isomorphism. Then for any $E \in \Fun(\Lambda,k)$, the map
$j^{o*}_l \to j^{o*} \circ \pi_{l!}$ adjoint to the isomorphism
$j^{o*}_l \circ \pi_l^* \cong j^{o*}$ induces a natural map
$$
HH_\idot(E) \cong H_\idot(\Delta^o,j^{o*}i_l^*E) \to
H_\idot(\Delta^o,j^{o*}\pi_{l!}i_l^*E) = HH_\idot(\pi_{l!}i_l^*E).
$$
Taking $E = W_nA^\hush$ and $l=p$, and composing this map with the
Verschiebung map \eqref{V.la}, we obtain a functorial map
\begin{equation}\label{V.hw}
V:W_nHH_\idot(A) \to HH_\idot(\pi_{p!}i_p^*W_nA^\hush) \to
W_{n+1}HH_\idot(A),
\end{equation}
the Hochschild-Witt homology version of the Verschiebung map.

\subsection{Comparison.}\label{cmp.subs}

The first of the comparison theorems for Hochschild-Witt homology
proved in \cite{ka.hw} works for any associative unital
$\F_p$-algebra $A$, but only in degree $0$. Here is the statement.

\begin{theorem}[{{\cite[Theorem 5.4]{ka.hw}}}]\label{hw.1.thm}
For any $n \geq 1$ and unital associative $\F_p$-algebra $A$, there
exists a functorial isomorphism $W_nHH_0(A) \cong W_n^H(A)$, where
$W_n^H(A)$ is the non-commutative Witt vectors groups introduced by
L. Hesselholt in \cite{hewi}. These isomorphisms are compatible with
the restriction maps $R$ of \eqref{R.hw} and Verschiebung maps
\eqref{V.hw}.\endproof
\end{theorem}

In particular, as a part of the construction, Hesselholt's Witt
vector groups fit into short exact sequences
\begin{equation}\label{he.seq}
\begin{CD}
A/[A,A] @>{V^n}>> W_{n+1}^H(A) @>{R}>> W_n^H(A) @>>> 0,
\end{CD}
\end{equation}
where $A/[A,A] = HH_0(A) = W_1^H(A)$ is the quotient of $A$ by the
subspace spanned by commutators. It has been further noted in
\cite{heerr} that the sequence is not in general exact on the
left. In terms of Theorem~\ref{hw.1.thm}, the sequence
\eqref{he.seq} is the degree-$0$ part of the homology long exact
sequence induced by \eqref{W.la.seq}, and the fact that it is not
exact on the left corresponds to the fact that the connecting
differential does not have to be trivial (with an explicit example
given in \cite{heerr}).

The second comparison theorem works in all degrees but only for
commutative algebras. To state it, recall that for any smooth
algebraic variety $X$ over $\F_p$, the {\em de Rham-Witt complex}
$W_\idot\Omega_X^\hdot$ of sheaves on $X$ has been constructed in
\cite{ill}. This is actually a projective system of complexes
$W_n\Omega_X^\hdot$, $n \geq 1$ of sheaves on $X$ in the Zariski
topology, functorial with respect to $X$, and equipped with
restriction maps $R:W_{n+1}\Omega_X^\hdot \to
W_n\Omega^\hdot_X$. The first term $W_1\Omega_X^\hdot$ is the usual
de Rham complex, and for any $n$, we have $W_n\Omega^0_X \cong
W_n(\calo_X)$, the ring of Witt vectors of the structure sheaf
$\calo_X$. The complexes $W_\idot\Omega^\hdot_X$ are also equipped
with functorial Verschiebung maps $V:W_n\Omega_X^\hdot \to
W_{n+1}\Omega_X^\hdot$ that extend the usual Verschibung maps on
$W_\idot\Omega^0_X$. We note that $V$ does not commute with the
differential --- instead, one has $pdV=Vd$.

\begin{theorem}[{{\cite[Theorem 6.14]{ka.hw}}}]\label{hw.2.thm}
Let $k=\F_p$, and let $A$ be an $\F_p$-algebra satisfying the
assumptions of Theorem~\ref{hkr.thm}. Then we have functorial
isomorphisms
\begin{equation}\label{hkr.hw}
W_nHH_i(A) \cong H^0(X,W_n\Omega^i_X), \qquad n \geq 1, i \geq 0
\end{equation}
that are compatible with the maps $R$, $V$, and send the
Connes-Tsygan differential $B$ to the de Rham-Witt differential
$d$.\endproof
\end{theorem}

We note that for a commutative $A$, Hesselholt's Witt vectors
$W_n^H(A)$ coincide with the classical Witt vectors. The
identification \eqref{hkr.hw} for $i=0$ is then the same as in
Theorem~\ref{hw.1.thm} (and this is in fact used in the proof of
Theorem~\ref{hw.2.thm}). Theorem~\ref{hw.2.thm} is similar to
Theorem~\ref{hkr.thm} in that it establishes Hochschild-Witt
homology groups as the correct non-commutative generalization of the
de Rham-Witt forms. However, it also clarifies somewhat the
structure of the original de Rham-Witt complex of Bloch, Deligne and
Illusie. For example, the limit $W\Omega_X^\hdot$ of the projective
system $W_\idot\Omega_X^\hdot$ has a natural decreasing filtration,
so that for any $n \geq 2$, we have a short exact sequence of
complexes
$$
\begin{CD}
0 @>>> \gr^nW\Omega_X^\hdot @>>> W_n\Omega^\hdot_X @>>>
W_{n-1}\Omega^\hdot_X @>>> 0.
\end{CD}
$$
In degree $0$, $\gr^nW\Omega^0_X \cong \calo_X$ does not depend on
$n$, but in higher degrees, the sheaves $\gr^n\Omega^i_X$ are new
functorial sheaves on $X$ that cannot be expressed in the terms of
the cotangent bundle and its tensor powers. This is true already for
$n=2$. Namely, denote by $\B^\hdot_X \subset \Zz^\hdot_X \subset
\Omega^\hdot_X$ the subsheaves of locally exact resp.\ locally
closed forms, and recall that we have the functorial Cartier
isomorphism
\begin{equation}\label{C.eq}
C:\Zz^\hdot_X/\B^\hdot_X \cong \Omega^\hdot_X.
\end{equation}
Then Illusie proves in \cite{ill} that one has a functorial short
exact sequence
\begin{equation}\label{gr.2}
\begin{CD}
0 @>>> \Omega^\hdot_X/\B^\hdot_X @>>> \gr^2W\Omega_X^\hdot @>>>
\B^\hdot_X @>>> 0.
\end{CD}
\end{equation}
Thus both $\Omega^\hdot_X = \gr^1W\Omega^\hdot_X$ and
$\gr^2W\Omega^\hdot_X$ have a two-step filtration with the same
associated graded quotients, but what is a subobject in one is the
quotient object in the other, and vice versa. The extension class
represented by the sequence \eqref{gr.2} is a new characteristic
class of smooth algebraic varieties over $\F_p$; there is no way to
construct it, short of doing the whole rather opaque construction of
the de Rham-Witt complex.

However, if one uses Theorem~\ref{hw.2.thm}, then the quotients
$\gr^nW\Omega_X^\hdot$ become much more transparent. Namely, a
non-commutative generalization of the Cartier isomorphism
\eqref{C.eq} has been given in \cite{ka.car} for odd $p$ and then in
\cite[Section 4]{ka.dege} for $p=2$. Using these results, one easily
shows that for $A$ as in Theorem~\ref{hw.2.thm}, the connecting
differential in the Hochschild homology long exact sequence induced
by \eqref{W.la.seq} vanishes, so that we have a natural isomorphism
$$
H^0(X,\gr^nW\Omega_X^i) \cong HH_i(\pi_{p^n!}i_{p^n}^*A^\hush).
$$
For $n=2$, the tilting-type relation between $\Omega_X^\hdot$ and
$\gr^2\Omega^\hdot_X$ is straightforwardly explained in terms of the
non-commutative Cartier map of \cite{ka.car}. For all $n$, the
behavior is similar; it is eludicated in \cite[Remark 6.13]{ka.hw}.

\section{Further developments and open questions.}\label{add.sec}

\subsection{Frobenius maps.}

Let us now describe some additional structures carried by Witt
vectors, both commutative and non-commutative. We start with the
Frobenius map. In the classical situation, the statement is as
follows.

\begin{prop}\label{F.prop}
There exist a unique collection of maps $F:W_{n+1}(A) \to W_n(A)$,
$n \geq 1$ that are additive, functorial with respect to the
commutative ring $A$, and satisfy $F \circ R = R \circ F$ and $w_n
\circ F = w_{n+1}$, $n \geq 1$. Moreover, these maps also satisfy $F
\circ V = p\id$ and $F(T(a)) = a^p$, $a \in A$, where $T:A \to
W_\idot(A)$ are the Teichm\"uller representative maps.
\end{prop}

\proof{} To find the maps $F$ such that $w_n \circ F = w_{n+1}$ and
$F \circ R = R \circ F$, it suffices to find universal polynomials
$f_i(a_0,\dots,a_i)$, $i \geq 1$ such that for any $n \geq 1$, we
have
\begin{equation}\label{F.eq}
f_1^{p^{n-1}} + pf_2^{p^{n-2}} + \dots + p^{n-1}f_n = a_0^{p^n} +
pa_1^{p^{n-1}} + \dots + p^na_n.
\end{equation}
To do this, assume by induction that we have found $f_i$ for $i <
n$. Then to show that \eqref{F.eq} provides a well-defined
polynomial $f_n$ with integer coefficients, it suffices to observe
that by Lemma~\ref{teich.bis.le}, $b^{p^{n-1}}=b^{p^n} \mod p^n$ for
any integer $n \geq 1$ and any integer $b$.

To prove that the resulting maps $F$ are unique, additive, and
satify all the equalities claimed, it suffices to consider the
universal situation $A=\Z$. But then the ghost maps $w_\idot =
\langle w_1,\dots,w_n \rangle:W_n(A) \to A^n$ are injective, so that
all the claims follow from the compatibility between $w_\idot$ and
$F$.
\endproof

We note that it also follows from \eqref{F.eq} that $f_n = a_{n-1}^p
\mod p$, $n \geq 1$. Therefore if $A$ consist entirely of
$p$-torsion, thus has a Frobenius endomorphism $\phi:A \to A$, $a
\mapsto a^p$, then $F:W_{n+1}(A) \to W_n(A)$ is given by $F =
W_n(\phi) \circ R$ and provides a canonical additive lifting of the
Frobenius map $\phi$ to Witt vectors. By abuse of terminology, $F$
is also called ``the Frobenius endomorphism'' for a general ring
$A$.

A special feature of the case $pA=0$ is the additional identity
$VF=p\id$. Indeed, this is an identity on universal polynomials with
coefficients in $\F_p$, so to check it, it suffices to check that it
holds after evaluation at elements in the algebraic closure
$\overline{\F}_p$. In other words, we have to prove that $VF=p\id$
for $A=\overline{\F}_p$. But then $\phi:\overline{\F}_p \to
\overline{\F}_p$ is a bijection, so that $F=W(\phi)$ is a bijection
as well, and $FVF=pF$ implies $VF=p\id$.

For polynomial Witt vectors $W_n(M)$, $M$ an $\F_p$-vector space,
defining a Frobenius map is even simpler. Namely, for any free
$Z$-module $N$ and integer $n \geq 1$, the tautological
identification $(N^{\otimes p^n}) \cong (N^{\otimes p})^{\otimes
  p^{n-1}}$ induces a map of functors $Q_n(N) \to Q_{n-1}(N^{\otimes
  p})$. This map is compatible with the trace functor structure, and
descends to a functorial map
\begin{equation}\label{F.W.eq}
F:W_n(M) \to W_{n-1}(M^{\otimes p})^{\overline{\sigma}}
\end{equation}
for any $\F_p$-vector space $M$, where $\overline{\sigma}$ is
induced by the trace functor structure. It is also immediately
obvious from construction that we have
$$
F \circ V = \id + \overline{\sigma} + \dots +
\overline{\sigma}^{p-1},
$$
where $V$ is the Verschiebung map \eqref{V.eq}. For any associative
unital $\F_p$-algebra $A$, the map \eqref{F.W.eq} induces a map
\begin{equation}\label{F.la}
F:W_nA^\hush \to \pi_{p*}i_p^*W_{n-1}A^\hush
\end{equation}
of functors from $\Lambda$ to abelian groups, and $F \circ
V:\pi_{p!}i_p^*W_nA^\hush \to \pi_{p*}i_p^*W_nA^\hush$ is the trace
map $\tr_p$ of \eqref{tr.l}. On the level of Hochschild-Witt
homology, we have maps
$$
F:W_{n+1}HH_\idot(A) \to W_nHH_\idot(A),
$$
and it has been proved in \cite[Corollary 1.10]{ka.hw} that
$\overline{\sigma}$ acts trivially on individual Hochschild homology
groups, so that $F \circ V = p\id:W_nHH_\idot(A) \to
W_nHH_\idot(A)$. Moreover, it has been proved in \cite[Lemma
  2.7]{ka.hw} that we have
\begin{equation}\label{fbv.eq}
FBV=B,
\end{equation}
where $B:W_nHH_\idot(A) \to W_nHH_{\idot+1}(A)$ is the Connes-Tsygan
differential.

As far as comparisons theorem of Subsection~\ref{cmp.subs} are
concerned, we note that Hesselholt does not directly construct a
Frobenius endomorphism on his Witt vector groups, but they inherit
it from the Topological Hochschild Homology spectrum, since he does
prove that the limit group $W^H(A)$ coincides with $\pi_0(TR(A,p))$
of \cite{BHM}. The identification of Theorem~\ref{hw.1.thm} is
probably compatible with the Frobenius maps but this has not been
checked. The de Rham-Witt complex $W\Omega_X^\hdot$ of a smooth
algebraic variety $X$ over $\F_p$ also has a Frobenius endomorphism
$F$; in fact, this is an integral part of the construction in
\cite{ill}. The endomorphisms $F$ and $V$ of the de Rham-Witt
complex commute, one has $FV=VF=p\id$, and moreover, one has
$FdV=d$, where $d$ is the differential in the complex. Our
identification of Theorem~\ref{hw.2.thm} is compatible with $F$; the
equality $FdV=d$ corresponds to the general non-commutative equality
\eqref{fbv.eq}.

\subsection{Multiplication.}

Next, let us describe the structure that we have omitted in
Theorem~\ref{witt.thm} --- namely, the structure of a commutative
ring.

\begin{theorem}\label{mult.thm}
There exists a unique collection of functorial unital commutative
ring structures on abelian groups $W_n(A)$, $n \geq 1$ such that the
restriction maps $R$ and the ghost maps $w_n$ are ring
maps. Moreover, with respect to these ring structures, we also have
$F(a \cdot b)=F(a) \cdot F(b)$, $V(a) \cdot b = V(a \cdot F(b))$,
$a,b \in W_n(A)$, and $T(ab) = T(a) \cdot T(b)$, $a,b \in A$.
\end{theorem}

\proof{} By induction, we may assume that we have constructed the
ring structures on $W_l(A)$, $l \leq n$. Decompose $W_{n+1}(A) =
A^{n+1} = A \times A^n \cong A \times W_n(A)$ so that
$R^n:W_{n+1}(A) \to A$ is the projection onto the first factor, and
define the product by
\begin{equation}\label{mult.eq}
\langle a_0,b_0 \rangle \cdot \langle a_1,b_1 \rangle = \langle
a_0a_1, b_0 \cdot F(\langle a_1,b_1\rangle) + a_0^pb_1 \rangle,
\end{equation}
where $F:W_{n+1}(A) \to W_n(A)$ is the Frobenius map of
Proposition~\ref{F.prop}. Then it immediately follows from the
inductive assumption that $w_n$ is multiplicative with respect to
the product \eqref{mult.eq}. Now as in Proposition~\ref{F.prop}, to
check that \eqref{mult.eq} defines a commutative associative unital
product compatible with the additive structure and satisfying all
the identities, it suffices to consider the case $A=\Z$, where
everything follows from compatibility with the ghost map. The same
goes for the uniqueness claim.
\endproof

We note that as immediate corollary of Theorem~\ref{mult.thm}, we
see that the image $VW(A) \subset W(A)$ of the Verschiebung map is
an ideal in the ring $W(A)$, and we have $A = W(A)/VW(A)$. If
$pA=0$, then $VF=p\id$ implies that $pW(A) \subset W(A)$ is
contained in $VW(A)$. If moreover $A$ is perfect --- that is, the
Frobenius map $\phi:A \to A$ is a bijection --- then $VW(A)=pW(A)$, so
that $W(A)$ is a $p$-adic lifting of the ring $A$ in the most naive
sense.

For polynomial Witt vectors, the product survives in the form of
external multiplication. Namely, for any two free $\Z$-modules
$N_0$, $N_1$ and integer $l$, we have a canonical isomorphism $(N_1
\otimes N_2)^{\otimes l} \cong N_1^{\otimes l} \otimes N_2^{\otimes
  l}$, and this induces a functorial map
\begin{equation}\label{mu.eq}
\mu:W_n(M_0) \otimes W_n(M_1) \to W_n(M_0 \otimes M_1)
\end{equation}
for any $\F_p$-vector spaces $M_0$, $M_1$ and integer $n \geq 1$.

\begin{prop}\label{mu.prop}
The map \eqref{mu.eq} is commutative, associative and unital in the
obvious sense. Moreover, we have
\begin{equation}\label{mu.RFV}
R \circ \mu = \mu \circ (R \otimes R), \ F \circ \mu = \mu
\circ (F \otimes F), \ \mu \circ (V \otimes \Id) = V \circ \mu
\circ (\Id \circ F).
\end{equation}
\end{prop}

\proof{} See \cite[Lemma 1.5, Proposition 3.10]{ka.witt}. \endproof

\begin{corr}
For any $\F_p$-vector space $M$ and integer $n$, we have
$$
V \circ F = p\id:W_n(M) \to W_n(M).
$$
\end{corr}

\proof{} Take $M_0=M$, $M_1=\F_p$, and apply \eqref{mu.RFV} and the
identity $VF=p\Id$ in $W(\F_p)=\Z_p$.
\endproof

Here is another application of the product \eqref{mu.eq}. Assume
given a finite-dimensional $\F_p$-vector space $M$, and let $M^*$ be
the dual vector space. Then for any $n \geq 1$, the natural pairing
$M \otimes M^* \to \F_p$ combined with the product map \eqref{mu.eq}
induces a functorial map
$$
W_n(M) \otimes W_n(M^*) \to W_n(\F_p) \cong \Z/p^n\Z.
$$
One shows (see \cite[Lemma 3.12]{ka.witt}) that this is in fact a
perfect pairing, and $W_n(M)$ and $W_n(M^*)$ are dual
$\Z/p^n\Z$-modules (note that typically they are not flat, but since
$\Z/p^n\Z$ is Gorenstein, duality is well-behaved). The duality
interchanges the Verschiebung map $V$ and the Frobenius map $F$. In
fact, it has been shown in \cite[Lemma 3.7]{ka.witt} that the short
exact sequence \eqref{W.seq} has a functorial dual sequence
$$
\begin{CD}
0 @>>> W_n(M) @>{C^m}>> W_{m+n}(M) @>{F^n}>> W_m(M^{\otimes p^n})^\tau
@>>> 0,
\end{CD}
$$
where $C:W_n(M) \to W_{n+1}(M)$ is a certain functorial {\em
  corestriction map} dual to the restriction map $R$ and satisfying
$RC=CR=p\id$.

In the simplest non-trivial case $n=2$, the picture can be explained
as follows. For any $\F_p$-vector space $M$, one has a natural
functorial exact sequence
\begin{equation}\label{4.term}
\begin{CD}
0 @>>> M @>{\psi}>> (M^{\otimes p})_\sigma @>{\tr_p}>> (M^{\otimes
  p})^\sigma @>{\wh{\psi}}>> M @>>> 0,
\end{CD}
\end{equation}
where $\psi$ sends $m \in M$ to $m^{\otimes p}$, and $\wh{\psi}$ is
dual. This sequence represents a certain extension class in the
category of functors from $\F_p$-vector spaces to $\F_p$-vector space,
and this class is non-trivial, so that \eqref{4.term} does not admit
a functorial splitting (for more details on this, see \cite[Section
  6]{ka.bok}). However, the extension becomes trivial if we consider
it in the category of functors from $\F_p$-vector spaces to abelian
groups, and $W_2(M)$ is exactly the splitting object. It has a
functorial three-step filtration $F^\hdot$ with associated graded
terms $M$, $\Phi(M)$, $M$, where $\Phi(M)$ is the image of the map
$\tr_p$ in \eqref{4.term}, and with terms $F^2W_1(M) \cong M$,
$F^1W_2(M) \cong (M^{\otimes p})_\sigma$, $F^0W_2(M)=W_2(M)$, while
the quotients are given by $W_2(M)/F^1W_2(M) \cong M$,
$W_2(M)/F^2W_2(M) \cong (M^{\otimes p})^\sigma$. In other words,
$W_2(M)$ can represented as either an extension of $M$ by
$(M^{\otimes p})_\sigma$, or as an extension of $(M^{\otimes
  p})^\sigma$ by $M$. As an abelian group, $W_2(M)$ is annihilated
by $p^2$ but not by $p$; multiplication by $p$ acts by the
composition
$$
\begin{CD}
W_2(M) @>{R}>> W_1(M) \cong M @>{C}>> W_2(M)
\end{CD}
$$
of the projection $R$ onto the top quotient $F^0/F^1$ of the
filtration, and the embedding $C$ of the bottom piece $F^2$. For
higher $n$, the picture is similar, although the filtration now has
$n(n+1)/2$ associated graded pieces; for details, see \cite[Section
  3.1]{ka.witt}.

The map \eqref{mu.eq} is also compatible with the trace functor
structure. The precise meaning of this can be found in \cite[Section
  4.3]{ka.witt}, and the implication is that we have natural maps
$$
\mu:W_nA^\hush \otimes W_nB^\hush \to W_n(A \otimes B)^\hush
$$
for any two associative unital $\F_p$-algebras $A$, $B$.  Combining
these maps with the K\"unneth isomorphism, we obtain maps
\begin{equation}\label{mu.hh}
\mu:W_nHH_\idot(A) \lotimes W_nHH_\idot(B) \to W_nHH_\idot(A \otimes
B),
\end{equation}
where $W_nHH_\idot(A)$, $W_nHH_\idot(B)$ are tacitly understood as
objects in the derived category of abelian groups, and $\lotimes$ is
the derived tensor product. These maps also satisfy
\eqref{mu.RFV}. If $A$ is commutative, then the multiplication map
$A \otimes A \to A$ is an algebra map, so that $W_nHH_\idot(A)$
becomes a commutative associative algebra. In the situation of
Theorem~\ref{hw.2.thm}, this algebra structure is identified with
the standard algebra structure on the de Rham-Witt complex under the
isomorphisms \eqref{hkr.hw}.

\subsection{Extending the definition.}

So far, for simplicity, we have only defined polynomial Witt vectors
for vector spaces over $\F_p$. However, when defining $W_n(M)$, it
is not necessary to represent $M$ as a quotient $N/pN$ of a free
$\Z$-module $N$ --- it has been shown in \cite[Lemma 2.1]{ka.witt}
that already for a free $\Z/p^n\Z$-module $N$ with $M\cong N/pN$, we
have $W_n(M) \cong Q_n(N)$. Therefore if we fix a perfect field $k$
of characteristic $p$, we can consider the category of free modules
over $W_n(k)$, the ring of $n$-truncated Witt vectors of the field
$k$, and repeat Definition~\ref{Q.def} {\em verbatim} for such free
modules. Then Corollary~\ref{W.corr} still holds, and we obtain a
functor $W_n$ from $k$-vector spaces to abelian groups. All the
other results about polynomial Witt vectors also hold --- in fact,
\cite{ka.witt} and \cite{ka.hw} work in this larger generality from
the very beginning. In particular, we have the Verschiebung maps
$V$, the Frobenius maps $F$, the restriction maps $R$ and the
correstriction maps $C$ (although one has to keep in mind that these
maps are not necessarily $W(k)$-linear --- for example, $F$ is
semilinear with respect to the Frobenius endomorphism of $W(k)$). We
also have the Hochschild-Witt homology groups and the product
\eqref{mu.eq}. Theorem~\ref{hw.2.thm} also holds in this larger
generality, with $W_\idot\Omega^\hdot_X$ being the relative de
Rham-Witt complex over $k$.

Another and more substantial extension concerns algebraic varieties
that are not necessarily affine. One could try to incorporate this
case by considering sheaves of non-commutative algebras, but this is
not what arises in practical applications. In fact, it is a
well-established principle of non-commutative geometry that ``every
non-commutative variety is affine in the derived sense'', or in
other words, that the correct basic object to consider is an
associative unital DG algebra, or a small DG category (see
e.g.\ \cite{kel} or \cite[Section 2]{orlo}, \cite[Section 1]{orlo1}
for a detailed discussion of this point). What one would like to
have, then, is the theory of Hochschild-Witt homology groups
$W_\idot HH_\idot(A_\idot)$ of a DG algebra $A_\idot$ over $k$.

A machine providing such a theory has in fact been constructed in
\cite{ka.trace}. As an input, it needs a trace functor $F$ from
$k$-vector spaces to abelian groups that has two additional
properties: it should be {\em balanced} in the sense of
\cite[Definition 3.9]{ka.trace} and {\em localizing} in the sense of
\cite[Definition 5.5]{ka.trace}. Both propeties are closed under
extensions (\cite[Lemma 3.11]{ka.trace} and \cite[Lemma
  5.6]{ka.trace}) and hold for the cyclic power functor $M \mapsto
M^{\otimes n}_\sigma$, $n \geq 1$ (\cite[Lemma 5.7, Proposition
  5.10]{ka.trace}). Therefore by \eqref{W.seq}, they also hold for
$W_n$, $n \geq 1$. What the machine produces is a collection of
twisted Hochschild homology functors $FHH_\idot(A_\idot)$ for small
DG categories $A_\idot$ over $k$. These functors are derived-Morita
invariant in the sense of \cite[Section 4.6]{kel} and moreover, are
additive invariants in the sense of \cite[Section
  5.1]{kel}. Plugging in the polynomial Witt vectors functors $W_n$,
we obtain Hochschild-Witt homology functors $W_nHH_\idot(A_\idot)$,
$n \geq 1$, that are additive invariants, just as they should
be. They are also related by the maps $F$, $V$, $R$, $C$, and carry
the external product \eqref{mu.hh}.

Unfortunately, the machine of \cite{ka.trace} is rather indirect and
cumbersome (it works by replacing a chain complex by a
simplicial-cosimplicial abelian group). There might be a more
straightforward approach to the problem that starts by extending
Definition~\ref{Q.def} directly to chain complexes. This has not
been done yet, although a way to do it for $W_2$ was sketched in
\cite[Section 6]{ka.bok}.

One further result that is not automatic at all is a K\"unneth-type
isomorphism for Hochschild-Witt homology. The difficulty here is
that it certainly does not hold on the level of Hochschild homology
groups --- nor de Rham-Witt forms, in the commutative case. For de
Rham-Witt forms, we do have K\"unneth isomorphism once we turn on the
differential (this gives cristalline cohomology). This suggests that
for general DG algebras $A_\idot$, $B_\idot$, the map \eqref{mu.hh}
should induce an isomorphism on periodic cyclic homology groups. At
present, this has not been checked.

Another interesting question that comes up in applications to DG
algebras is the degeneration of the Hodge-to-de Rham spectral
sequence \eqref{hdr.bis.seq}. For DG algebras $A_\idot$ over a ring
$k$ that contains $\Q$, it has been conjectured by Kontsevich and
Soibelman that \eqref{hdr.bis.seq} degenerates at first term if
$A_\idot$ is homologically smooth and homologically proper (these
properties correspond to the usual smoothness and properties of
algebraic varieties, see \cite{kel}, \cite[Section 2]{orlo} for more
details). The conjecture has been proved recently in \cite{ka.dege},
and an even more general statement has been discovered in
\cite{efi}. The method of the proof is that of Deligne and Illusie
\cite{DI}, and it works essentially by reducing the question to
positive characteristic. However, the argument of \cite{DI} in fact
appeared as a simplification of a degeneration for de Rham-Witt
complex established already in \cite{ill}. Thus it is natural to
expect that the sequence \eqref{hdr.bis.seq} for Hochschild-Witt
cohomology also degenerates when the DG category $A_\idot$ is smooth
and proper --- either on the nose, or at least, when one takes the
limit with respect to $n$ and inverts $p$. However, at present, this
has not been studied at all.

\subsection{Big Witt vectors.}

Yet another generalization of the Witt vctors definition is actually
quite old, and exists already on the classical level --- this the
theory of ``big'' or ``universal'' Witt vectors that combines
together the theories for all primes $p$.

Classically, for any unital commutative ring $A$, one consider the
product $\W(A) = A^{\N}$ of copies of $A$ numbered by positivie
integers $m \geq 1$, and defines the universal ghost map
$\wh{w}_\idot:\W(A) \to A^{\N}$ by
\begin{equation}\label{gh.b.eq}
\wh{w}_m(a_1,\dots,a_m) = \sum_{d | m}da_d^{m/d},
\end{equation}
where the sum is over all the integers $d \geq 1$ that divide
$m$. We note that for any prime $p$ and integer $n \geq 1$, the map
$\wh{w}_{p^n}$ only depends on the components $a_{p^i}$, $0 \leq i
\leq n$, and up to a renumbering of these components, it coincides
with the ghost map of \eqref{gh.eq}, called {\em $p$-typical} in
this context. One then equips $\W(A)$ with the product topology and
proves that there exists a unique functorial continuous commutative
ring structure on $\W(A)$ such that for any $m \geq 1$,
\eqref{gh.b.eq} is a ring map.

The proof can be done along the same line as Theorem~\ref{witt.thm},
but there is an attractive alternative: the additive group structure
on $\W(A)$ is visible right away. Namely, one shows that there is a
functorial isomorphism
\begin{equation}\label{w.t}
\W(A) \cong 1+tA[[t]] \subset A[[t]]^*,
\end{equation}
where $A[[t]]^*$ is the group of invertible formal power series in
one formal variable $t$ with coefficients in $A$, considered with
respect to multiplication. Explicitly, the isomorphism
\eqref{w.t} sends an element $\langle a_1,a_2,\dots \rangle \in
\W(A)$ to the series
$$
\sum_{i \geq 1}(1-a_it^i) \in 1+tA[[t]].
$$
The fact that this is an isomorphism is obvious by induction, and
the fact that \eqref{gh.b.eq} is additive is a simple computation.

The product in the ring $\W(A)$ is not immediately obvious from the
isomorphism \eqref{w.t}. To see it, one can use an interpretation
due to Almkvist \cite{A} that exhibits a dense subgroup in $\W(A)$
spanned by characteristic polynomials of matrices over $A$, with
addition corresponding to the direct sum of matrices -- and then
shows that the product corresponds to the tensor
product. Alternatively, it has been shown in \cite{ka.jap} how to
write down the product on the whole $\W(A)$ in terms of the Tate
residue (or more precisely, the Contou-Carr\`ere residue, see
e.g.\ \cite{go1}, \cite{go2}, \cite{go3} and references therein). To
obtain this interpretation, one starts by observing that $\W(A)$ is
the kernel of the split surjective map
$$
K_1(A[[t]]) \to K_1(A)
$$
of the first algebraic $K$-groups of $A[[t]]$ and $A$. This also
suggests a way to obtain a polynomial version of the functor
$\W$. Namely, for any perfect field $k$ and $k$-vector space $M$,
let $\wh{T}^\hdot(M/k)$ be the completion of the tensor algebra
$T^\hdot(M/k)$ generated by $M$ over $k$ with respect to the
augmentation ideal $T^{\geq 1}(M/k) \subset T^\hdot(M/k)$. Then we
also have a split surjective map
$$
K_1(\wh{T}^\hdot(M/k)) \to K_1(k),
$$
and one can define $\W(M/k)$ as the kernel of this map. It seems
that this gives a good object equipped with all the additional
structures one would like, such as Frobenius and Verschiebung maps
for each prime $p$, a product similar to \eqref{mu.eq}, and a trace
functor structure that allows one to define the universal
Hochschild-Witt homology groups. However, at present, there is no
written proof of all this, and we will return to the topic in the
forthcoming paper \cite{ka.polywitt}.

To connect the universal and the $p$-typical theory, one can compute
the universal Witt vectors ring $\W(\Z)$. As a group, we have
$$
\W(\Z) = \Z\langle \eps_1,\eps_2, \dots \rangle,
$$
the group of possibly infinite linear combinations of elements
$\eps_n$, $n \geq 1$. The product is given by
$$
\eps_i\eps_j = \frac{ij}{\{i,j\}}\eps_{\{i,j\}},
$$
where $\{i,j\}$ is the least common multiple of the integers $i,j
\geq 1$. Then by functoriality, $\eps_i$ acts on $\W(A)$ for any
commutative ring $A$, and if $A$ is $p$-local --- that is, any
integer $n$ prime to $p$ acts on $A$ by an invertible map --- then
$(1/n)\eps_n$ is an idempotent endomorphism of the ring
$\W(A)$. Taken together, these idempotents generate a
decomposition
$$
\W(A) = \prod_nW(A)
$$
into the product of copies of $W(A)$ numbered by integers $n \geq 1$
prime to $p$ (this is known as the {\em $p$-typical
  decomposition}). For polynomial Witt vectors, the same
decomposition seems to exist; to construct it, one needs to use the
technology of $\Z$-Mackey profunctors developed in \cite[Section
  9.2]{ka.ma}.

\noindent
{\sc Steklov Mathematics Institute,\\
Moscow, Russia}

\medskip

\noindent
{\em E-mail address\/}: {\tt kaledin@mi.ras.ru}

\end{document}